\documentclass[a4paper,reqno,12pt,psamsfonts]{amsart}

\usepackage[headings]{fullpage}

\usepackage{upref}

\usepackage{hyperref}

\usepackage{amsthm}
\usepackage{amsmath}
\usepackage{amssymb}

\usepackage[all]{xy}   

\theoremstyle{plain}
\newtheorem{thm}{Theorem}[section]
\newtheorem*{thm*}{Theorem}
\newtheorem{prop}[thm]{Proposition}
\newtheorem{lem}[thm]{Lemma}
\newtheorem{cor}[thm]{Corollary}

\theoremstyle{definition}
\newtheorem{defn}[thm]{Definition}

\theoremstyle{remark}
\newtheorem*{rem*}{Remark}
\newtheorem*{notation*}{Notation}
\newtheorem*{ex*}{Example}
\newtheorem*{ackn*}{Acknowledgements}

\theoremstyle{definition}
\newtheorem*{comp_ax*}{Comparison axiom}
\newtheorem*{ext_ax*}{Extension axiom}
\newtheorem*{weak_comp_ax*}{Weak comparison axiom}
\newtheorem*{cov_ax*}{Covering axiom}
\newtheorem*{zero_ax*}{Zero axiom}
\newtheorem*{str_ext_ax*}{Strong extension axiom}

\newcommand{\Id}{\operatorname{Id}}
\newcommand{\Vol}{\operatorname{Vol}}
\newcommand{\sys}{\operatorname{sys}}

\begin{document}

\title[Homological invariance]{Homological invariance for asymptotic invariants and systolic inequalities}
\author{Michael Brunnbauer}
\address{Mathematisches Institut, Ludwig-Maximilians-Universit\"at M\"unchen, Theresienstr.\ 39, D-80333 M\"unchen, Germany}
\email{michael.brunnbauer@mathematik.uni-muenchen.de}
\date{\today} 
\keywords{systolic constant, minimal volume entropy, spherical volume}
\subjclass[2000]{Primary 53C23;   
                 Secondary 53C20, 
                           57R19} 

\begin{abstract}
We show that the systolic constant, the minimal volume entropy, and the spherical volume of a manifold depend only on the image of the fundamental class under the classifying map of the universal covering. Moreover, we compute the systolic constant of manifolds with fundamental group of order two (modulo the value for the real projective space) and derive an inequality between the minimal volume entropy and the systolic constant.
\end{abstract}

\maketitle



\section{Introduction}\label{intro}

In the present article, we will consider connected closed smooth manifolds $M$. We will prove that certain asymptotic and systolic invariants of $M$ depend only on the image of the fundamental class under the classifying map of the universal covering. 

These invariants include the \emph{minimal volume entropy} $\lambda(M)$, which describes the asymptotic volume growth of the universal covering, the \emph{spherical volume} $T(M)$, which is an invariant intermediate between the minimal volume entropy and the simplicial volume, and the \emph{systolic constant} $\sigma(M)$, which determines the relation between the lengths of short non-contractible loops and the volume of the manifold. We will show the following theorem, which is the main result of this paper:

\begin{thm}\label{main_thm_intro}
Let $M$ and $N$ be two connected closed smooth manifolds having the same fundamental group $\pi$. Let $\Phi:M\to K(\pi,1)$ and $\Psi:N\to K(\pi,1)$ denote the respective classifying maps of the universal coverings of $M$ and $N$. If the subgroups of orientation preserving loops of $M$ and $N$ coincide and if $\Phi_*[M] = \Psi_*[N]$, then 
\[ I(M)=I(N), \] 
where $I$ denotes either the systolic constant $\sigma$, the minimal volume entropy $\lambda$, or the spherical volume $T$.
\end{thm}

Here, the fundamental classes $[M]$ and $[N]$ are to be understood with respect to (local) coefficients in the orientation bundles of $M$ and $N$. Note that in the oriented case the condition on the orientation preserving subgroups is always fulfilled.

Most cases of this theorem for the systolic constant were known by work of Babenko (see \cite{Babenko(2006)}), whose ideas we follow in parts of the proof. Sabourau applied these ideas to the minimal entropy (see \cite{Sabourau(2006)}). The spherical volume has not been considered in this context up to now.

To unify the treatment of the mentioned invariants, we will introduce certain axioms that are satisfied by $\sigma$, $\lambda$, and $T$. In the proof of Theorem \ref{main_thm_intro} we will use only these axioms and no other properties of the invariants. Thus the theorem holds for all invariants $I$ fulfilling the axioms. One more example for such an invariant is the \emph{stable systolic constant} $\sigma^{st}$, a variation of the systolic constant. Moreover, it will be convenient to consider relative versions of the invariants (relative to some homomorphism $\phi:\pi_1(M)\to\pi$ from the fundamental group to an arbitrary group) and to extend their definitions to simplicial complexes. (The definitions and the axioms can be found in sections \ref{sys+entr}, \ref{axioms}, and \ref{spher_vol}. In section \ref{systolic_mfd} a homological criterion for the vanishing of the systolic constant is given. Section \ref{main_proof} contains the proof of Theorem \ref{main_thm_intro}.) 

As an application of Theorem \ref{main_thm_intro}, we will look at manifolds whose fundamental groups have only two elements. The fact that $K(\mathbb{Z}_2,1)=\mathbb{R}P^\infty$ will allow us to derive a complete list of possible values for the systolic constant. (For the minimal volume entropy and the spherical volume this case is of no interest since these invariants vanish for finite fundamental groups.)

\begin{cor}\label{fund_grp_two_elem_intro}
Let $M$ be a connected closed $n$-dimensional smooth manifold. If $\pi_1(M)=\mathbb{Z}_2$, then
\[ \sigma(M) = 
\begin{cases}
\sigma(\mathbb{R}P^n) & \Phi_*[M]_{\mathbb{Z}_2}\neq 0 \\
0 & \mbox{otherwise}
\end{cases}, \]
where $\Phi:M\to\mathbb{R}P^\infty$ denotes the classifying map of the universal covering of $M$.
\end{cor}

This was previously known only for orientable manifolds by another paper of Babenko (see \cite{Babenko(2004)}). Note also that the exact value of $\sigma(\mathbb{R}P^n)$ is unknown except for dimension two, where it is $2/\pi$ (see \cite{Pu(1952)}). This corollary will be proved in section \ref{appl}.

In section \ref{str_ext_ax_section} we will investigate what happens to the minimal volume entropy when one enlarges the fundamental group by attaching $1$-cells to the manifold. Using these observations, the main theorem, and the computations from \cite{Sabourau(2006)}, we will finally prove in the last section:

\begin{thm}\label{sabourau_ineq_intro}
Let $M$ be a connected closed $n$-dimensional smooth manifold. Then there exists a positive constant $c_n$ depending only on $n$ such that
\[ \sigma(M) \geq c_n \frac{\lambda(M)^n}{\log^n(1+\lambda(M))}. \]
\end{thm}

The proof will cover the dimensions $n\geq 3$. For surfaces, the theorem was shown in \cite{KS(2005)} by Katz and Sabourau. Moreover, Sabourau proved this inequality in special cases including the case of aspherical orientable manifolds (see \cite{Sabourau(2006)}). Note that the proof of Theorem \ref{sabourau_ineq_intro} requires relative versions of both invariants involved in the formulation.

Theorem \ref{sabourau_ineq_intro} sharpens a theorem of Gromov (\cite{Gromov(1983)}, Theorem 6.4.D') that stated the inequality 
\[ \sigma(M) \geq c'_n \frac{\|M\|}{\log^n(1+\|M\|)} \]
for oriented manifolds, where $\|M\|$ denotes the \emph{simplicial volume} of $M$. Namely, recall that there is another inequality by Gromov (\cite{Gromov(1982)}, pages 35-37), improved by \cite{BCG(1991)}, Th\'eor\`emes 3.8 and 3.16, that says
\[ {\textstyle\frac{n^{n/2}}{n!}}\|M\| \leq \lambda(M)^n. \]
So indeed, Theorem \ref{sabourau_ineq_intro} implies Gromov's inequality (up to constants).

Note also that these two lower bounds for the systolic constant are optimal in dimensions two and three. A discussion of this optimality result may be found in \cite{KSV(2007)}.

Most of the results in this text hold also for piecewise linear manifolds. But for simplicity and readability, we will only consider smooth manifolds. Since our interest lies in invariants that were originally defined for smooth manifolds, we think this restriction is justifiable. So henceforth, all manifolds in this text will be supposed smooth.

The next section contains the technical core of this article. It is concerned with maps from manifolds to CW complexes that will be deformed by elementary methods to gain useful normal forms for such maps.

\begin{ackn*}
I would like to thank D.\,Kotschick for his continuous advice, help, and encouragement. I am also grateful to the referee for useful remarks and valuable hints. Financial support from the \emph{Deutsche Forschungsgemeinschaft} is gratefully acknowledged.
\end{ackn*}


\section{Topological preliminaries}\label{top_prelim}

In this section we show that any map from a manifold to a CW complex can be brought into a form convenient for our purposes. In the last subsection we define the notion of absolute degree and point out its geometrical meaning.

\subsection{The Hopf trick}

Consider a proper map $f:(M,\partial M)\to (X,A)$ from a manifold $M$ (with or without boundary) to a relative manifold $(X,A)$, both having dimension $n$. Recall that a pair $(X,A)$ is called \emph{relative manifold} if $X$ is a Hausdorff space and $A\subset X$ is a closed subspace such that $X\setminus A$ is a manifold (see \cite{Spanier(1966)}, page 297). For example, every $n$-dimensional CW complex is a manifold relative to its $(n-1)$-skeleton.

Let $p\in X\setminus A$ be a point. Replacing $f$ by a properly homotopic map if necessary, we may assume that $f$ is smooth on the preimage of a small neighborhood of $p$, and moreover that $f$ is transverse to $p$. The preimage of $p$ then consists of finitely many points $p_1,\ldots,p_\ell$ in $M\setminus\partial M$. Choosing a local orientation of $X$ at $p$, the map $f$ induces local orientations of $M$ at these points.

In this situation the following `trick' due to Hopf applies, see \cite{Hopf(1930)}. A modern presentation can be found in \cite{Epstein(1966)}, pages 378-380. (There $X$ is supposed to be a manifold. But in fact, it is enough that $X$ is a manifold around $p\in X$.)

\begin{lem}[Hopf trick]
Let $n\geq 3$. Assume now, that there is a path $\gamma$ in $M$ between two preimage points, say from $p_1$ to $p_2$, that reverses the induced orientations and that is mapped to a contractible loop in $X$. Then we may deform $f$ on a compact subset of $M\setminus\partial M$ such that the number $\ell$ of preimage points of $p$ is reduced by $2$.
\end{lem}

\subsection{Orientation issues}

If $M$ is a connected compact manifold of dimension $n$, then there is exactly one local integer coefficient system $\mathcal{O}$ on $M$ such that $H_n(M,\partial M;\mathcal{O})\cong\mathbb{Z}$, namely the orientation bundle of $M$. It is trivial ($\mathcal{O}=\mathbb{Z}$) if and only if $M$ is orientable. For all other $\mathbb{Z}$ coefficient systems $H_n(M,\partial M;\mathcal{O})=0$.

Moreover note that each local integer coefficient system on a connected CW complex $X$ is determined by a unique homomorphism 
\[ \rho:\pi_1(X)\to\mathbb{Z}_2=\operatorname{Aut}(\mathbb{Z}). \]
Denote this coefficient system by $\mathcal{O}_\rho$. In the case of the orientation bundle of a manifold, $\rho$ has kernel the subgroup of orientation preserving loops. If not said otherwise, we will always use this homomorphism $\rho:\pi_1(M)\to\mathbb{Z}_2$.

\begin{rem*}
With respect to coefficients in the orientation bundle $\mathcal{O}_\rho$ all paths in $M$ are orientation preserving.
\end{rem*}

\subsection{Maps to $n$-dimensional CW complexes}

Consider a map $f:(M,\partial M)\to (X,A)$ to a pair of CW complexes whose induced homomorphism $f_*$ on fundamental groups is surjective. If $\ker f_*\subset\ker\rho$, then the homomorphism $\rho:\pi_1(M)\to\mathbb{Z}_2$ induces a homomorphism $\rho:\pi_1(X)\to\mathbb{Z}_2$ and the homomorphisms on homology
\[ f_*: H_*(M,\partial M;\mathcal{O}_\rho) \to H_*(X,A;\mathcal{O}_\rho) \]
are well-defined. If it is possible, we will always use coefficients in $\mathbb{K}=\mathcal{O}_\rho$, but if $\ker f_*\not\subset\ker\rho$, we have to take $\mathbb{K}=\mathbb{Z}_2$ coefficients. 

\begin{rem*}
Note that $\ker f_*\subset \ker\rho$ if and only if the covering $\tilde M_{f_*}$ associated to the subgroup $\ker f_*\subset\pi_1(M)$ is orientable.
\end{rem*}

Assume for the rest of this subsection that $X$ is $n$-dimensional and $A\subset X$ is an $(n-1)$-dimensional subcomplex. Then
\[ H_n (X,A;\mathbb{K}) \cong \ker (H_n(X,X^{(n-1)};\mathbb{K})\xrightarrow{\partial} H_{n-1}(X^{(n-1)},A;\mathbb{K})) \]
by the long exact homology sequence of the triple $(X,X^{(n-1)},A)$. Moreover by excision $H_n(X,X^{(n-1)};\mathbb{K})$ is isomorphic to 
\[ \bigoplus_{\mbox{\scriptsize $e$ $n$-cell}} \mathbb{Z}\cdot e, \;\;\;\;\mbox{respectively}\;\; \bigoplus_{\mbox{\scriptsize $e$ $n$-cell}} \mathbb{Z}_2\cdot e. \]

Let $a\in H_n(X,A;\mathbb{K})$ be given by $\sum_{\mbox{\scriptsize $e$ $n$-cell}} r_e\cdot e$ with $r_e= 0$ for all but finitely many $n$-cells $e$.

\begin{lem}\label{main_top_lemma}
Let $n\geq 3$. If $f:(M,\partial M)\to (X,A)$ fulfills $f_*[M]_\mathbb{K} =a$ and is surjective on fundamental groups, then it is homotopic relative to $\partial M$ to a map $f':(M,\partial M)\to(X,A)$ such that there are $\sum_{\mbox{\scriptsize\rm{$e$ $n$-cell}}} |r_e|$ many pairwise disjoint closed balls $D_{ei_e},i_e=1,\ldots,|r_e|,$ in $M\setminus\partial M$ such that 
\[ f'^{-1}(\mathring{e})=\mathring{D}_{e1}\cup\ldots\cup\mathring{D}_{e|r_e|} \] 
and 
\[ f':\mathring{D}_{e1}\cup\ldots\cup\mathring{D}_{e|r_e|}\to\mathring{e} \] 
is a covering of $\mathbb{K}$ degree $r_e$, i.\,e.\ it is an $r_e$-sheeted covering such that the $\mathbb{K}$ orientations on the balls $\mathring{D}_{ei}$ agree.
\end{lem}

\begin{notation*}
The absolute value on $\mathbb{Z}_2$ is defined as $0$ for the trivial and $1$ for the non-trivial element.
\end{notation*}

\begin{proof}
First, remove the interiors of all $n$-cells $e$ of $X$ with $f(M)\cap\mathring{e}=\emptyset$ and with $r_e=0$. This affects neither the surjectivity of $f_*:\pi_1(M)\twoheadrightarrow\pi_1(X)$ nor the equality $f_*[M]_\mathbb{K}=a$. (If $X'$ denotes the complex obtained from $X$ by removing those open $n$-cells, then $H_n(X',A;\mathbb{K}) \hookrightarrow H_n(X,A;\mathbb{K})$ is injective by the long exact sequence of the tripel $(X,X',A)$ and the second statement follows.) Thus, there remain only finitely many $n$-cells because $M$ is compact and $r_e=0$ for almost all $n$-cells $e$.

We proceed by induction on the number of remaining $n$-cells of $X$. If it is zero, there is nothing to prove. 

Now let $e$ be one of the $n$-cells of $X$. Choose a point $p\in\mathring{e}$ and assume w.\,l.\,o.\,g.\ that $f$ is transverse to it. Denote its preimages by $p_1,\ldots,p_\ell$. The assumption implies that $f$ has local $\mathbb{K}$ degree $r_e$ at $p$. Hence $\ell\geq |r_e|$.

In case $\ker f_*\subset\ker\rho$, we may choose $d:=\ell-|r_e|$ points from the points $p_1,\ldots,p_\ell$ such that one half of them is mapped orientation preservingly to $p$, the other half orientation reversingly (with respect to some choice of orientation of $X$ at $p$). Take a path $\alpha$ from a point of the first kind to one of the second kind. Since $f_*:\pi_1(M)\twoheadrightarrow\pi_1(X)$ is surjective the loop $f\circ\alpha$ lies in its image. Let $\beta$ be a loop based at the first point that is mapped to $f\circ\alpha$. Then $\gamma := \beta^{-1}\alpha$ has contractible image under $f$ and is orientation reversing with respect to any choice of local orientations of $M$ at $p_1,\ldots,p_\ell$ coming from $p$. Hence we may apply the Hopf trick and reduce $d$ by two. By induction $d$ finally becomes zero and $\ell=|r_e|$. 

Now consider the case that $\ker f_*\not\subset\ker\rho$. Choose a local orientation for $X$ at $p$ and thus also for $M$ at the $p_i$. Let $\alpha$ be a path between two such points. Proceeding as above, we may assume that its image is contractible in $X$. If $\alpha$ is not orientation reversing, choose a loop $\beta$ based at the starting point of $\alpha$ that is orientation reversing and in the kernel of $f_*$. Then $\gamma:=\beta\alpha$ is orientation reversing and is mapped to a contractible loop. The Hopf trick reduces $\ell$ by two and induction shows that we may deform $f$ until $\ell=|r_e|$.

Next, choose a ball $D\subset \mathring{e}$ with center $p$ such that $f^{-1}(D)$ consists of $|r_e|$ pairwise disjoint closed balls $D_{e1},\ldots,D_{e|r_e|}$ in $M\setminus\partial M$, each mapped diffeomorphically and with the same orientation behaviour onto $D$. Compose $f$ with a strong deformation retraction from $X\setminus\mathring{D}$ to $X\setminus\mathring{e}$.

The induction hypothesis applied to $M':=M\setminus (\mathring{D}_{e1}\cup\ldots\cup\mathring{D}_{e|r_e|})$, $X' := X\setminus\mathring{e}$ and $a' := \sum_{e'\neq e}r_{e'}\cdot e'$ finishes the proof.
\end{proof}

\subsection{Maps to arbitrary CW complexes}

Let $M$ be a connected closed manifold of dimension $n\geq 2$. Consider a map $f:M\to X$ to a CW complex $X$ that is surjective on fundamental group level. Let $a\in H_n(X^{(n)};\mathbb{K})$ be a homology class (where we use $\mathbb{K}=\mathcal{O}_\rho$ if $\ker f_*\subset\ker\rho$ and $\mathbb{K}=\mathbb{Z}_2$ otherwise as always) and let $i:X^{(n)}\hookrightarrow X$ be the inclusion.

\begin{lem}\label{n-skeleton}
If $f_*[M]_\mathbb{K} = i_* a$, then we may deform $f$ such that its image lies in the $n$-skeleton of $X$ and that $f_*[M]_\mathbb{K}=a\in H_n(X^{(n)};\mathbb{K})$.
\end{lem}

In the case $\mathbb{K}=\mathcal{O}_\rho$ this is due to Babenko, see \cite{Babenko(2006)}, Lemme 3.10.

\begin{proof}
First consider the case $\mathbb{K}=\mathcal{O}_\rho$. The Hurewicz theorem gives an epimorphism
\[ h:\pi_{n+1}(X^{(n+1)},X^{(n)}) \twoheadrightarrow H_{n+1}(X^{(n+1)},X^{(n)};\mathcal{O}_\rho) \]
since $\pi_k(X^{(n+1)},X^{(n)})=0$ for $k\leq n$. Hence the commutative diagram with exact rows and vertical Hurewicz homomorphisms
\[\xymatrix{
\pi_{n+1}(X^{(n+1)},X^{(n)}) \ar@{->>}[d]_-h \ar[r]^-\partial & \pi_n(X^{(n)}) \ar[d]_-h \ar[r] & \pi_n(X^{(n+1)})=\pi_n(X) \\
H_{n+1}(X^{(n+1)},X^{(n)};\mathcal{O}_\rho) \ar[r]^-\partial & H_n(X^{(n)};\mathcal{O}_\rho) \ar[r]^-{j_*} & H_n(X^{(n+1)};\mathcal{O}_\rho) = H_n(X;\mathcal{O}_\rho)
}\]
shows that the kernel of $j_*$ is contained in the image of $h:\pi_n(X^{(n)})\to H_n(X^{(n)};\mathcal{O}_\rho)$. 

For the case $\mathbb{K}=\mathbb{Z}_2$ note that we may add a further row to the above diagram by applying the reduction map from $\mathcal{O}_\rho$ to $\mathbb{Z}_2$ coefficients. The induced map
\[ H_{n+1}(X^{(n+1)},X^{(n)};\mathcal{O}_\rho) \twoheadrightarrow H_{n+1}(X^{(n+1)},X^{(n)};\mathbb{Z}_2) \]
is obviously surjective. Hence we get the commutative diagram
\[\xymatrix@C-0pt{
\pi_{n+1}(X^{(n+1)},X^{(n)}) \ar@{->>}[d]_-h \ar[r]^-\partial & \pi_n(X^{(n)}) \ar[d]_-h \ar[r] & \pi_n(X) \\
H_{n+1}(X^{(n+1)},X^{(n)};\mathbb{Z}_2) \ar[r]^-\partial & H_n(X^{(n)};\mathbb{Z}_2) \ar[r]^-{j_*} & H_n(X;\mathbb{Z}_2)
}\]
and again we have $\ker j_* \subset \operatorname{im} h$.

By cellular approximation we may assume that $f$ maps to the $n$-skeleton of $X$. Since $f_*[M]_\mathbb{K}=i_*a\in H_n(X;\mathbb{K})$, we see that $j_*f_*[M]_\mathbb{K} = j_*a$, hence
\[ a - f_*[M]_\mathbb{K} \in \ker j_*. \]
Let $s:S^n\to X^{(n)}$ be a preimage under $h$, i.\,e.\ we have $s_*[S^n]_\mathbb{K}=a- f_*[M]_\mathbb{K}$. We may assume that $s$ is contractible in $X$. In fact, we can choose $s$ in the image of the boundary homomorphism $\partial:\pi_{n+1}(X^{(n+1)},X^{(n)}) \to \pi_n(X^{(n)})$. Define
\[ f': M\to M\vee S^n \xrightarrow{f\vee s} X^{(n)}, \]
where the first map contracts the boundary of a small ball in $M$. Then $f'_*[M]_\mathbb{K} = f_*[M]_\mathbb{K} +s_*[S^n]_\mathbb{K} = a\in H_n(X^{(n)};\mathbb{K})$ and the maps $f$ and $f'$ are homotopic as maps to $X$ by the choice of $s$.
\end{proof}

\subsection{Absolute and geometric degree}

Let $f:(M,\partial M)\to (N,\partial N)$ be a map between two connected compact manifolds of dimension $n$. It factors as
\[ (M,\partial M)\xrightarrow{\bar f} (\bar N,\partial \bar N) \xrightarrow{p} (N,\partial N) \]
where $p:\bar N\to N$ is the covering map associated to the subgroup $f_*(\pi_1(M))\subset \pi_1(N)$. Let $j$ be the number of sheets of $p$. If $\ker \bar f_*=\ker f_*\subset \ker\rho$, then we may define the \emph{degree} of $f$ as zero for $j=\infty$ and as $j\cdot\deg(\bar f)$ for $j<\infty$ where $\deg(\bar f)$ is determined by
\[ \bar f_* [M]_{\mathcal{O}_\rho} = \deg(\bar f)\cdot [\bar N]_{\mathcal{O}_\rho}. \]
This is to be understood as $\deg(\bar f)=0$ if $H_n(\bar N,\partial \bar N;\mathcal{O}_\rho)=0$. (This degree is defined only up to sign. We have to choose orientations to get a well-defined integer.)

Moreover, we define the \emph{absolute degree} of $f$ by
\[ \deg_a(f) := 
\begin{cases}
0 & j=\infty \\
j\cdot |\!\deg(\bar f)| & j<\infty,\ker f_*\subset \ker\rho \\
j\cdot |\!\deg_2(\bar f)| & j<\infty,\ker f_*\not\subset \ker\rho
\end{cases},\]
where $\deg_2(\bar f)$ denotes the $\mathbb{Z}_2$ degree of $\bar f$. (This number is well-defined without any choices.)

\begin{rem*}
This definition coincides with the usual definition of absolute degree (see e.\,g.\ \cite{Epstein(1966)} or \cite{Skora(1987)}).
\end{rem*}

The \emph{geometric degree} $\deg_g(f)$ of $f$ is the smallest integer $d$ for which there is a map $f':(M,\partial M)\to (N,\partial N)$ homotopic to $f$ that is transverse to some point $p\in N\setminus\partial N$ such that $f'^{-1}(p)$ consists of $d$ points. Note that always $\deg_a(f)\leq \deg_g(f)$.

\begin{thm*}[Hopf, Kneser]
If $n\geq 3$, then $\deg_a(f) = \deg_g(f)$. In the two dimensional case the same equality holds if one assumes that $M$ and $N$ are closed.
\end{thm*}

\begin{proof}[Proof for $n\geq 3$]
Choose a CW decomposition of $(N,\partial N)$ and lift it to $(\bar N,\partial\bar N)$. If the number of sheets of $p$ is infinite, then $\bar N$ is not compact and therefore $H_n(\bar N,\partial\bar N;\mathbb{K})=0$. Hence $\bar f_*[M]_\mathbb{K}=0$ and Lemma \ref{main_top_lemma} shows that $\bar f$, and consequently $f$, contracts to the respective $(n-1)$-skeleton. In particular $\deg_g(f)=0$.

Now assume $j<\infty$. Applying Lemma \ref{main_top_lemma} to $\bar f$ with $a=\deg_a(\bar f)\cdot [\bar N]_\mathbb{K}$, we get a homotopic map $\bar f'$ such that each open $n$-cell of $\bar N$ is covered by exactly $\deg_a(\bar f)$ open $n$-cells in $M$. Hence $f':=p\circ \bar f'$ is homotopic to $f$ and has geometric degree equal to $j\cdot\deg_a(\bar f')=\deg_a(f')$.
\end{proof}

The Hopf part ($n\geq 3$) of this theorem was proved in \cite{Hopf(1930)} by using the Hopf trick, see also \cite{Epstein(1966)}. In \cite{Epstein(1966)}, it is stated incorrectly that the equivalence of absolute and geometric degree also holds without further assumptions for $n=2$. See \cite{Skora(1987)} for a discussion of this and a modern proof of Kneser's result from \cite{Kneser(1930)}.

Our proof in fact shows more, namely that each topdimensional open cell of $N$ is covered by exactly $\deg_a(f)$ open cells in $M$. If we use the fact that smooth manifolds are triangulable, we get maps having the following nice property:

\begin{defn}
A simplicial map $f:X\to Y$ between two $n$-dimensional simplicial complexes is said to be \emph{$(n,d)$-monotone} if the preimage of every open $n$-simplex of $Y$ consists of at most $d$ open $n$-simplices in $X$. It is called \emph{strictly $(n,d)$-monotone} if the preimage of every open $n$-simplex of $Y$ consists of exactly $d$ open $n$-simplices in $X$.
\end{defn}

\begin{rem*}
Usually a map is called \emph{monotone} if the preimage of any point is connected. In this sense, $(n,1)$-monotone means that the map $f:X\to Y$ is monotone outside the $(n-1)$-skeleton and $(n,d)$-monotone means that we may divide $X\setminus X^{(n-1)}$ into $d$ sets such that $f$ is monotone on each of these sets.
\end{rem*}

\begin{cor}\label{deg_mon}
Let $f:(M,\partial M) \to (N,\partial N)$ be a map between connected compact manifolds of dimension $n\geq 3$ and let $d:=\deg_a(f)$. Then $f$ is homotopic to a strictly $(n,d)$-monotone map. In the two dimensional case one has to assume that $M$ and $N$ are closed to get the same conclusion.
\end{cor}

In the closed case this corollary was proved without any dimensional restrictions in \cite{Babenko(1992)}, proof of Proposition 2.2, part (a). The two-dimensional case of this corollary is proved by using Kneser's theorem. In fact, Kneser constructed in his original proof a strictly $(2,d)$-monotone map homotopic to $f$.


\section{Systolic constants and minimal volume entropy}\label{sys+entr}

In this section we introduce relative versions of the systolic constant and of the minimal volume entropy. Moreover, we define a stable version of the systolic constant.

Let $X$ be a connected finite simplicial complex of dimension $n$ and let $\phi:\pi_1(X)\to\pi$ be a group homomorphism. There is a corresponding map $\Phi:X\to K(\pi,1)$ that induces this homomorphism $\phi$ on fundamental groups and that is determined uniquely up to homotopy by this property. By a Riemannian metric on $X$ we always understand a continuous piecewise smooth one.

\begin{defn}
For a Riemannian metric $g$ on $X$ define the \emph{$\phi$-systole $\sys_\phi(X,g)$} as the infimum of all lengths of closed piecewise smooth curves in $X$ whose images under the corresponding map $\Phi:X\to K(\pi,1)$ are non-contractible. The \emph{systolic constant relative to $\phi$} is given by
\[ \sigma_\phi(X) := \inf_g \frac{\Vol(X,g)}{\sys_\phi(X,g)^n}, \]
where the infimum is taken over all Riemannian metrics $g$ on $X$.
\end{defn}

\begin{rem*}
The systolic constant $\sigma_\phi(X)$ depends only on the kernel of $\phi$, not on $\phi$ itself. But for practical reasons we will keep $\phi$ in the notation. Nevertheless note that it is no actual restriction to assume that $\phi$ is surjective. The same remark applies to the next two definitions and to the definition of the spherical volume (Definition \ref{def_spher_vol}).
\end{rem*}

\begin{defn}
Denote by $H(\pi)$ the free Abelian group
\[ H_1(\pi;\mathbb{Z})/\mbox{Torsion}. \]
For an element $a\in H(\pi)$ denote by $\ell_g(a)$ the length of the shortest loop $\gamma$ in $X$ such that $\Phi'\circ\gamma\simeq a$ in $K(H(\pi),1)$ where $\Phi':X\to K(H(\pi),1)$ corresponds to the composition of $\phi$ and the canonical epimorphism $\pi\twoheadrightarrow H(\pi)$. Setting
\[ |a|_g := \lim_{k\to\infty} \frac{\ell_g(ka)}{k}, \]
we define the \emph{stable $\phi$-systole} by
\[ \sys^{st}_\phi(X,g) := \inf_{a\neq 0} |a|_g. \]
The \emph{stable systolic constant relative to $\phi$} is given by
\[ \sigma^{st}_\phi(X) := \inf_g \frac{\Vol(X,g)}{\sys^{st}_\phi(X,g)^n}. \]
\end{defn}

\begin{defn}
In the setting above let $\tilde X_\phi$ be the Galois covering of $X$ associated to the normal subgroup $\ker\phi\lhd\pi_1(X)$. For any Riemannian metric $g$ on $X$ define the \emph{volume entropy relative to $\phi$} as
\[ \lambda_\phi(X,g) := \lim_{r\to\infty} {\textstyle\frac{1}{r}} \log\Vol(B(x,r)), \]
where $B(x,r)$ is the ball of radius $r$ around a point $x\in\tilde X_\phi$ with respect to the lifted metric. This limit exists and is independent of the center $x\in\tilde X_\phi$ (see \cite{Manning(1979)}). One defines the \emph{minimal volume entropy of $X$ relative to $\phi$} as
\[ \lambda_\phi (X) := \inf_g \lambda_\phi (X,g)\Vol(X,g)^{1/n}. \]
\end{defn}

To work with this definition quickly gets complicated. But there is an equivalent definition that is easier to handle:

\begin{defn}\label{norms}
Let $G$ be a finitely generated group. A \emph{norm} on $G$ is a non-negative function $L:G \to [0,\infty)$ such that
\begin{enumerate}
\item $L(g)=0 \;\;\Leftrightarrow\;\; g=1$.
\item $L(g^{-1}) = L(g)$.
\item $L(gg') \leq L(g) + L(g')$ (triangle inequality).
\end{enumerate}
The \emph{growth function} $\beta_L : [0,\infty)\to [1,\infty]$ of a norm $L$ is defined by
\[ \beta_L(r) := \# \{g\in G | L(g)\leq r \}. \]
If the limit
\[ \lambda(G,L) := \lim_{r\to\infty}{\textstyle\frac{1}{r}}\log\beta_L(r) \]
exists, it is called the \emph{volume entropy} of $G$ with respect to $L$.
\end{defn}

\begin{rem*}
We may use the inclusion $\iota_x:\Gamma\hookrightarrow \tilde X_\phi,\gamma\mapsto\gamma\cdot x$ of the Galois group $\Gamma:=\pi_1(X)/\ker\phi$ into the Galois covering to induce a norm $L_{g,x}$ on $\Gamma$. Then it can easily be seen that
\[ \lambda_\phi (X,g) = \lambda(\Gamma,L_{g,x}) \]
by using translates of a fundamental domain of the Galois action. (See e.\,g.\ \cite{KH(1995)}, Proposition 9.6.6 or \cite{Sabourau(2006)}, Lemma 2.3.)
\end{rem*}


\section{Axioms for invariants}\label{axioms}

To unify the treatment of the systolic constants and the minimal volume entropy, we will investigate real-valued invariants $I_\phi(X)$ that are defined for connected finite simplicial complexes $X$ and group homomorphisms $\phi:\pi_1(X)\to\pi$ and that fulfill certain axioms. The main examples are $I=\sigma$ and $I=\lambda^n$.

\subsection{Comparison axiom and homotopy invariance}\label{comp_axiom}

For manifolds, invariants fulfilling the following comparison axiom behave reasonably well with respect to the absolute degree and are invariant under homotopy equivalence.

\begin{comp_ax*} 
Let $X$ and $Y$ be two connected finite simplicial complexes of dimension $n$ and let $\phi:\pi_1(X)\to\pi$ and $\psi:\pi_1(Y)\to\pi$ be group homomorphisms. If there exists an $(n,d)$-monotone map $f:X\to Y$ such that $\phi=\psi\circ f_*$, then
\[ I_\phi(X) \leq d\cdot I_\psi(Y). \]
\end{comp_ax*}

\begin{lem}[Babenko, Sabourau]\label{comp_ax_lambda_sigma}
The comparison axiom is fulfilled by $I=\sigma$, $I=\sigma^{st}$, and $I=\lambda^n$ (i.\,e.\ the minimal volume entropy to the power of $n$).
\end{lem}

Proofs of this lemma may be found in \cite{Babenko(2006)}, Proposition 3.2 and \cite{Sabourau(2006)}, Lemma 3.5 (both for $d=1$) and also in \cite{Babenko(1992)}, Propositions 2.2 and 8.7 (where $f_*:\pi_1(X)\to\pi_1(Y)$ is assumed to be surjective).

From Corollary \ref{deg_mon} we deduce:

\begin{cor}\label{homotopy_invar}
Let $M$ and $N$ be two connected closed manifolds and let $\psi:\pi_1(N)\to\pi$ be a group homomorphism. Let $f:M\to N$ be a map and $d:=\deg_a(f)$. If $I$ fulfills the comparison axiom, then
\[ I_{\psi\circ f_*}(M) \leq d\cdot I_\psi(N). \]
In particular, if $f$ is a homotopy equivalence, then
\[ I_{\psi\circ f_*} (M) = I_\psi (N). \]
\end{cor}

\subsection{Extension axiom}\label{ext_ax_section}

Since, given two manifolds, there may be no non-trivial map (say of absolute degree one) between them, we need a procedure to enlarge one of them such that we get a non-trivial map from the other manifold to the enlarged one. The next axiom shows how to do this and what happens to the invariants in the process.

Let $h:S^{k-1}\to X$ be a simplicial map, $1\leq k<n$, such that $\Phi\circ h$ is contractible if $k=2$. Define $X':=X\cup_h D^k$. This can be considered as a simplicial complex such that $X$ is a subcomplex. Define moreover $\phi':\pi_1(X)\to\pi$ as $\phi$ for $k\geq 3$ (the fundamental group has not changed), as the quotient map for $k=2$, and as an arbitrary extension of $\phi$ for $k=1$. Then we have $\phi'\circ i_* = \phi$ where $i:X\hookrightarrow X'$ is the inclusion.

\begin{defn}
An \emph{extension} $(X',\phi')$ of $(X,\phi)$ is a simplicial complex that is obtained by a finite number of attachments in the way described above.
\end{defn}

\begin{ext_ax*}
Let $(X',\phi')$ be an extension of $(X,\phi)$ where $\phi:\pi_1(X)\twoheadrightarrow\pi$ is an epimorphism. Then
\[ I_{\phi'}(X') = I_\phi(X). \]
\end{ext_ax*}

\begin{rem*}
The surjectivity assumption on $\phi$ guarantees that the corresponding Galois group $\pi\cong\pi_1(X)/\ker\phi$ remains the same for every extension of $(X,\phi)$. (Otherwise it could become bigger by attaching $1$-cells, see section \ref{str_ext_ax_section}.) In other words, the covering $\tilde X'_{\phi'}$ is obtained from $\tilde X_\phi$ by $\pi$-equivariantly attaching cells of dimension $1\leq k<n$.
\end{rem*}

\begin{lem}[\cite{Babenko(2006)}, Proposition 3.6 and \cite{Sabourau(2006)}, Lemma 3.6]\label{ext_ax_lambda_sigma}
The extension axiom is fulfilled in the cases $I=\sigma$, $I=\sigma^{st}$, and $I=\lambda$.
\end{lem}

Here Babenko and Sabourau restrict to $1<k<n$. But since $\phi$ is surjective, we only have to make the attached $1$-cells long enough until they play no role for the (stable) $\phi$-systole (longer than $\sys_\phi(X,g)$, respectively than $\sys^{st}_\phi(X,g)$) or for the norm $L_{g,x}$ (longer than some based loop in $X$ that represents the same element in $\pi$).


\section{Systolic manifolds}\label{systolic_mfd}

Using the topological lemmata from section \ref{top_prelim} and a famous theorem by Gromov, we are now able to give a homological classification of \emph{$\phi$-systolic} manifolds, i.\,e.\ of those manifolds $M$ with $\sigma_\phi(M)>0$.

\begin{defn}
A connected finite $n$-dimensional simplicial complex $X$ is called \emph{geometrically $\phi$-essen\-tial} for a group homomorphism $\phi:\pi_1(X)\to\pi$ if the associated map $\Phi:X\to K(\pi,1)$ does not contract to the $(n-1)$-skeleton of $K(\pi,1)$.
\end{defn}

In his Filling paper Gromov proved the following systolic inequality (\cite{Gromov(1983)}, Appendix 2, (B'$_1$)):

\begin{thm*}[Gromov]
If $X$ is geometrically $\phi$-essential, then
\[ \sys_\phi(X,g) \leq C_n\cdot\Vol(X,g)^{1/n} \]
for some universal constant $C_n>0$.
\end{thm*}

So geometrical $\phi$-essentiality implies $\phi$-systolicity. Using the comparison axiom we see immediately that this is an equivalence:
\[ \sigma_\phi(X) >0 \;\;\;\Leftrightarrow\;\;\; \mbox{$X$ geometrically $\phi$-essential}. \]

From the Lemmata \ref{n-skeleton} and \ref{main_top_lemma} follows directly:

\begin{cor}
A connected closed manifold $M$ of dimension $n\geq 3$ is not geometrically $\phi$-essential if and only if $\Phi_*[M]_\mathbb{K}=0$ in $H_n(\pi;\mathbb{K})$, where $\phi:\pi_1(X)\twoheadrightarrow\pi$ is assumed to be surjective.
\end{cor}

Together with the above equivalence this implies the following classification of $\phi$-systolic manifolds:

\begin{cor}\label{systolic_classification}
Let $M$ be a connected closed manifold of dimension $n\geq 3$ and let $\phi:\pi_1(M)\twoheadrightarrow\pi$ be an epimorphism. Then
\[ \sigma_\phi(M)>0 \;\;\;\Leftrightarrow\;\;\;
\begin{cases}
\Phi_*[M]_{\mathcal{O}_\rho} \neq 0 & \tilde M_\phi\;\mbox{orientable} \\
\Phi_*[M]_{\mathbb{Z}_2} \neq 0 & \tilde M_\phi\;\mbox{non-orientable}
\end{cases}.\]
\end{cor}

In the case where $\tilde M_\phi$ is orientable, this corollary is due to Babenko (see \cite{Babenko(1992)}, Theorem 8.2).


\section{Spherical volume}\label{spher_vol}

In this section we want to investigate another invariant: the spherical volume $T$. Its definition is a bit more involved than the definitions of the minimal volume entropy and the systolic constants. Therefore it is not easy to prove that $T$ fulfills the comparison axiom of section \ref{comp_axiom}. But we are able to find weaker axioms that lead to the same conclusions and thus show that Corollary \ref{homotopy_invar} is also valid in the case $I=T$.

The original definition of the spherical volume is due to Besson, Courtois, and Gallot (\cite{BCG(1991)} and \cite{BCG(1995)}, see also \cite{Storm(2002)}). Inspired by the definition of the minimal volume entropy and the systolic constant we introduce a relative version.

Again, let $X$ be a connected finite simplicial complex of dimension $n$ and let $\phi:\pi_1(X)\to\pi$ be a group homomorphism. Denote by $\tilde X_\phi$ the Galois covering associated to the normal subgroup $\ker\phi\lhd\pi_1(X)$ and by $\Gamma:=\pi_1(X)/\ker\phi$ the Galois group.

\begin{defn}\label{def_spher_vol}
Let $g$ be a Riemannian metric on $X$. Then $L^2(\tilde X_\phi)$ denotes the Hilbert space of square-integrable functions on $\tilde X_\phi$ with respect to the Riemannian volume of the lifted metric and $S^\infty(\tilde X_\phi)\subset L^2(\tilde X_\phi)$ denotes its unit sphere. Note that $\Gamma$ acts isometrically on both spaces by $\gamma\cdot f(x) := f(\gamma^{-1}x)$. Let $\mathcal{N}$ consist of those maps $F:\tilde X_\phi\to S^\infty(\tilde X_\phi)$ that are $\Gamma$-equivariant, Lip\-schitz continuous, and non-negative, i.\,e.\ their values are non-negative functions. If $F\in\mathcal{N}$, then its restriction to the interior of the $n$-cells is differentiable almost everywhere by Rademacher's theorem and we can define
\[ g_x^F (v_1,v_2) := \langle D_x F(v_1), D_x F(v_2) \rangle_{L^2(\tilde X_\phi)} \]
for almost all $x\in\tilde X_\phi, v_1,v_2\in T_x\tilde X_\phi$. (Tangent spaces are well-defined for points inside topdimensional simplices.) One finds that $g^F$ is an almost everywhere defined positive semi-definite $\Gamma$-invariant metric on $\tilde X_\phi$.

This metric descends to $X$ where it is also called $g^F$. We may define its volume form as $0$ at points where $g^F$ is degenerate or not defined and as the usual volume form at points where it is non-degenerate. Then $dV_{g^F}$ is an integrable $n$-form on $X$. Hence we can define 
\begin{align*}
\Vol(X,g^F) &:= \int_X dV_{g^F} \;\;\;\mbox{and} \\
T_\phi(X) &:= \inf_{F\in\mathcal{N}} \Vol(X,g^F).
\end{align*}
This number is called the \emph{spherical volume of $X$ relative to $\phi$}.
\end{defn}

\begin{rem*}
This definition is independent of the choice of the Riemannian metric $g$ on $X$ since the Hilbert spaces $L^2(\tilde X_\phi)$ for different Riemannian metrics are $\Gamma$-equivariantly isometric. Moreover, the notion of Lipschitz continuity of $F:\tilde X_\phi\to S^\infty(\tilde X_\phi)$ also does not depend on which metric $g$ we choose because $X$ is compact.
\end{rem*}

\begin{weak_comp_ax*} 
Let $X$ and $Y$ be connected finite simplicial complexes of dimension $n$ and let $\phi:\pi_1(X)\to\pi$ and $\psi:\pi_1(Y)\to\pi$ be group homomorphisms. If there exists a strictly $(n,d)$-monotone map $f:X\to Y$ such that $\phi=\psi\circ f_*$ and such that the induced homomorphism $f_*:\pi_1(X)/\ker\phi\xrightarrow{\scriptscriptstyle\cong}\pi_1(Y)/\ker\psi$ between the Galois groups is an isomorphism, then
\[ I_\phi(X) \leq d\cdot I_\psi(Y). \]
\end{weak_comp_ax*}

\begin{lem}\label{w_comp_ax_T}
The weak comparison axiom is fulfilled for $I=T$.
\end{lem}

\begin{proof}
Let $g_2$ be a Riemannian metric on $Y$. Define a Riemannian metric $g_1$ on $X$ by using $f^*g_2$ on the non-degenerate simplices and extending it over all of $X$.

Since $f_*:\pi_1(X)/\ker\phi \xrightarrow{\scriptscriptstyle\cong} \pi_1(Y)/\ker\psi$ is an isomorphism, we get an equivariant lift $\tilde f:\tilde X_\phi \to \tilde Y_\psi$ of $f$ that is again strictly $(n,d)$-monotone. The map
\begin{align*} 
\mathcal{I}:L^2(\tilde Y_\psi) &\to L^2(\tilde X_\phi),\\ 
\varphi &\mapsto \tilde\chi/\sqrt{d} \cdot (\varphi\circ\tilde f),
\end{align*}
where $\tilde\chi:\tilde X_\phi\to\mathbb{R}$ is the characteristic map of the non-degenerate $n$-simplices, is an equivariant isometric homomorphism that preserves non-negativity.

If $F:\tilde Y_\psi\to S^\infty(\tilde Y_\psi)$ is non-negative equivariant Lipschitz, then consider the non-negative equivariant Lipschitz map
\[ \mathcal{I}\circ F\circ \tilde f: \tilde X_\phi \to S^\infty(\tilde X_\phi). \]
We have $g^{\mathcal{I}\circ F\circ\tilde f}=g^{F\circ\tilde f}$ since $\mathcal{I}$ is isometric and
\[ \Vol(X,g^{F\circ\tilde f}) = d\cdot \Vol(Y,g^F), \]
which can be seen by looking at each open $n$-simplex of $Y$ together with its preimage separately. Therefore
\[ T_\phi(X) \leq d\cdot T_\psi(Y). \qedhere \]
\end{proof}

\begin{cov_ax*} 
Let $f:X\to Y$ be a $d$-sheeted covering map of connected finite simplicial complexes and let $\psi:\pi_1(Y)\to\pi$ be a homomorphism. Then
\[ I_{\psi\circ f_*}(X) \leq d\cdot I_\psi(Y). \]
\end{cov_ax*}

\begin{lem}
The covering axiom is true for $I=T$.
\end{lem}

\begin{proof}
This proof is essentially the same as the proof of Lemma \ref{w_comp_ax_T} with one exception: the lifted map $\tilde f:\tilde X_{\psi\circ f_*}\to\tilde Y_\psi$ is a covering map with 
\[ D:=[\ker\psi : f_*(\ker(\psi\circ f_*))] \] 
sheets. Therefore we have to replace the factor $\tilde\chi/\sqrt{d}$ in the definition of the isometry $\mathcal{I}$ by $1/\sqrt{D}$. (There are no degenerate simplices.) Then everything works out well.
\end{proof}

\begin{zero_ax*} 
Let $X$ and $Y$ be two connected finite simplicial complexes of dimension $n$ and let $\psi:\pi_1(Y)\to\pi$ be a group homomorphism. If $f:X\to Y$ is $(n,0)$-monotone, then
\[ I_{\psi\circ f_*}(X) =0. \]
\end{zero_ax*}

\begin{lem}\label{T_zero_ax}
The zero axiom is valid for $I=T$.
\end{lem}

\begin{proof}
We have 
\[ 2^n n^{n/2} T_\phi(X) \leq \lambda_\phi (X)^n \] 
for all $(X,\phi)$. (See \cite{BCG(1991)}, Th\'eor\`eme 3.8 or \cite{Storm(2002)}, Proposition 4.1. There the inequality is neither stated for simplicial complexes nor in the relative case, but the proof is exactly the same.) 

Since in the setting of the zero axiom $\lambda_{\psi\circ f_*}(X)^n=0$ by the comparison axiom we get $T_{\psi\circ f_*}(X) =0$ from the cited inequality.
\end{proof}

Now we can prove that Corollary \ref{homotopy_invar} also holds for $I=T$:

\begin{prop}
Let $M$ and $N$ be two connected closed manifolds and let $\psi:\pi_1(N)\to\pi$ be a group homomorphism. Let $f:M\to N$ be a map and $d:=\deg_a(f)$. Then
\[ I_{\psi\circ f_*}(M) \leq d\cdot I_\psi(N) \]
for any invariant $I$ that fulfills the weak comparison axiom, the covering axiom, and the zero axiom. 
\end{prop}

\begin{proof}
Denote by $p:\bar N\to N$ the connected covering of $N$ associated to the subgroup $f_*(\pi_1(M))\subset \pi_1(N)$. If $\bar N$ is not compact, then $\deg_a(f)=0$ and $f$ is homotopic to a $(n,0)$-monotone map by Corollary \ref{deg_mon}. By the zero axiom $I_{\psi\circ f_*}(M)=0$.

Assume now that $\bar N$ is compact. Note that $f$ factorizes over $\bar N$
\[\xymatrix{
& \bar N \ar[d]^p \\
M \ar[r]_-f \ar[ru]^-{\bar f} & N
}\]
such that $\bar f_*:\pi_1(M)\twoheadrightarrow\pi_1(\bar N)$ is surjective and the absolute degree factors as $\deg_a(f) = \deg_a(\bar f)\deg_a(p)$. By Corollary \ref{deg_mon} we may homotope $\bar f$ to be strictly $(n,\deg_a(\bar f))$-mono\-tone. By the weak comparison axiom applied to this map and the covering axiom applied to $p$ the proposition follows.
\end{proof}

For homotopy invariance even less assumptions are needed:

\begin{cor}\label{homotopy_invar_2}
If $f:M\xrightarrow{\scriptscriptstyle\simeq} N$ is a homotopy eqivalence, then 
\[ I_{\psi\circ f_*}(M) = I_\psi(N) \]
for every invariant $I$ satisfying the weak comparison axiom.
\end{cor}

The extension axiom does not need to be adjusted:

\begin{lem}\label{ext_ax_T}
The invariant $I=T$ fulfills the extension axiom.
\end{lem}

\begin{proof}
We have $\pi_1(X)/\ker\phi \cong\pi\cong \pi_1(X')/\ker\phi'$, hence $\tilde X'_{\phi'}$ is obtained from $\tilde X_\phi$ by equivariant attachments of cells of dimension less than $n$. These cells are zero sets, thus canonically
\[ L^2(\tilde X'_{\phi'}) = L^2(\tilde X_\phi). \]
Restriction to $X$ defines a surjective map $\mathcal{N}'\to\mathcal{N}$. Since the non-negative part of $S^\infty(\tilde X_\phi)$ is contractible we may extend any map $F\in \mathcal{N}$ equivariantly over $\tilde X'_{\phi'}$ to get a map $F'\in\mathcal{N}'$. This gives a section $\mathcal{N}\to\mathcal{N}'$ of the surjection above. Furthermore $\Vol(X',g^{F'})=\Vol(X,g^F)$ because the attached cells are of lower dimension, hence zero sets. Thus $T_{\phi'}(X')=T_\phi(X)$.
\end{proof}


\section{Homological invariance}\label{main_proof}

Throughout this section let $I$ be an invariant that fulfills both the weak comparison axiom and the extension axiom. The main examples are $I\in\{\sigma, \lambda^n, T\}$. We now prove our main result, Theorem \ref{main_thm_intro}. In fact, we show an analogous statement for the relative case which includes Theorem \ref{main_thm_intro} as a special case.

\begin{thm}\label{main_thm}
Let $M$ and $N$ be two connected closed manifolds of dimension $n\geq 3$ and let $\phi:\pi_1(M)\twoheadrightarrow\pi$ and $\psi:\pi_1(N)\twoheadrightarrow\pi$ be two epimorphisms. Denote by $\Phi:M\to K(\pi,1)$ and $\Psi:N\to K(\pi,1)$ the associated maps. 
\begin{enumerate}
\item If there exists a homomorphism $\rho:\pi\to\mathbb{Z}_2$ such that both $H_n(M;\mathcal{O}_{\rho\phi})\neq 0$ and $H_n(N;\mathcal{O}_{\rho\psi})\neq 0$ and if moreover 
\[ \Phi_*[M]_{\mathcal{O}_{\rho\phi}} = \Psi_*[N]_{\mathcal{O}_{\rho\psi}} \in H_n(\pi;\mathcal{O}_\rho) \] 
holds, then $I_\phi(M)=I_\psi(N)$.
\item If $\tilde N_\psi$ is non-orientable and 
\[ \Phi_*[M]_{\mathbb{Z}_2} = \Psi_*[N]_{\mathbb{Z}_2} \in H_n(\pi;\mathbb{Z}_2), \] 
then $I_\phi(M)\geq I_\psi(N)$.
\end{enumerate}
\end{thm}

Note that part (i) of this theorem in the absolute case is exactly Theorem \ref{main_thm_intro} for $n\geq 3$. The two-dimensional case is trivial since two closed surfaces with the same fundamental group are diffeomorphic. For the systolic constants $\sigma$ and $\sigma^{st}$ most of (i) is known by work of Babenko (see \cite{Babenko(2006)}).

From case (ii) follows immediately: 

\begin{cor}
If both $\tilde M_\phi$ and $\tilde N_\psi$ are non-orientable and 
\[ \Phi_*[M]_{\mathbb{Z}_2} = \Psi_*[N]_{\mathbb{Z}_2} \in H_n(\pi;\mathbb{Z}_2), \]
then $I_\phi(M) = I_\psi(N)$.
\end{cor}

For future use (section \ref{ineq_sys_entr}) we will consider pseudomanifolds.

\begin{defn}[see \cite{Spanier(1966)}, page 148]
A \emph{connected closed $n$-dimensional pseudomani\-fold} $X$ is a finite simplicial complex such that every simplex is a face of an $n$-simplex, every $(n-1)$-simplex is the face of exactly two $n$-simplices and for every two $n$-simplices $s$ and $s'$ there exists a finite sequence $s=s_1,\ldots,s_m=s'$ of $n$-simplices such that $s_i$ and $s_{i+1}$ have an $(n-1)$-face in common.
\end{defn}

\begin{rem*}
Since $X$ admits a CW decomposition with exactly one $n$-cell (see \cite{Sabourau(2006)}, Lemma 2.2), we find that $H_n(X;\mathcal{O}_\rho)$ is either $0$ or isomorphic to $\mathbb{Z}$ depending on the homomorphism $\rho:\pi_1(X)\to\mathbb{Z}_2$. Since there is no notion of orientation preserving paths in $X$, there may be more than one homomorphism $\rho:\pi_1(X)\to\mathbb{Z}_2$ (or indeed none) with $H_n(X;\mathcal{O}_\rho)\cong\mathbb{Z}$. Nevertheless, $H_n(X;\mathbb{Z}_2)=\mathbb{Z}_2$ in any case.
\end{rem*}

To prove Theorem \ref{main_thm} we need the following topological theorem, whose proof uses almost everything of section \ref{top_prelim}:

\begin{thm}\label{exists_deg_1_map}
Let $X$ be a connected closed pseudomanifold of dimension $n\geq 3$ and let $N$ be a connected closed $n$-dimensional manifold. Let $\phi:\pi_1(X)\twoheadrightarrow\pi$ and $\psi:\pi_1(N)\twoheadrightarrow\pi$ be two epimorphisms and let $\Phi:X\to K(\pi,1)$ and $\Psi:N\to K(\pi,1)$ be the associated maps. 
\begin{enumerate}
\item If either there is a homomorphism $\rho:\pi\to\mathbb{Z}_2$ such that both $H_n(X;\mathcal{O}_{\rho\phi})\neq 0$ and $H_n(N;\mathcal{O}_{\rho\psi})\neq 0$ and 
\[ \Phi_*[X]_{\mathcal{O}_{\rho\phi}} = \Psi_*[N]_{\mathcal{O}_{\rho\psi}} \in H_n(\pi;\mathcal{O}_\rho), \] 
\item or if $\tilde N_\psi$ is non-orientable and 
\[ \Phi_*[X]_{\mathbb{Z}_2} = \Psi_*[N]_{\mathbb{Z}_2} \in H_n(\pi;\mathbb{Z}_2), \]
\end{enumerate}
then there exists an extension $(X',\phi')$ of $(X,\phi)$ and a strictly $(n,1)$-monotone map $h:N\to X'$ such that
\begin{gather*}
\psi = \phi'\circ h_* \;\;\mathit{and} \\
h_*[N]_\mathbb{K} = i_*[X]_\mathbb{K} \in H_n(X';\mathbb{K}),
\end{gather*}
where $i:X\hookrightarrow X'$ is the inclusion.
\end{thm}

\begin{proof}
Use $\phi$ to identify $\pi_1(X)/\ker\phi=\pi$. Choose a CW decomposition of $X$. Now attach (possibly infinitely many) $2$-cells to $X$ whose attaching loops generate $\ker\phi$. Thus we get a CW complex $X(2)$ that has fundamental group $\pi_1(X(2))=\pi$. Next attach $3$-cells to $X(2)$ and kill $\pi_2(X(2))$ and then $4$-cells and so on. We obtain a sequence $X\subset X(2)\subset X(3)\subset\ldots$ of CW complexes that fulfill
\[ \pi_1(X(k))=\pi \;\;\mbox{and}\;\; \pi_s(X(k))=0 \;\;\mbox{for}\;\; 2\leq s<k. \]
This gives a CW decomposition of $K(\pi,1)$ and we have
\begin{align*}
X(n-1) &= K(\pi,1)^{(n-1)} \cup \left( \bigcup_{\mbox{\scriptsize$e$ $n$-cell of $X$}} e \right) \;\;\;\mbox{and} \\
X(k) &= K(\pi,1)^{(k)} \;\;\;\;\;\mbox{for}\; k\geq n.
\end{align*}

By Lemma \ref{n-skeleton}, $\Psi$ gives a map
\[ g:N\to X(n) \]
such that $g_*[N]_\mathbb{K}=i_*[X]_\mathbb{K}$. 

Lemma \ref{main_top_lemma} shows that we may deform $g$ to
\[ \hat g:N\to X(n-1) \]
with $\hat g_*[N]_\mathbb{K} = i_* [X]_\mathbb{K}$. Moreover, $\hat g$ is strictly $(n,1)$-monotone.

By compactness, we may choose a finite subcomplex $X\subset X'\subset X(n-1)$ such that $\hat g(N)\subset X'$ and $\hat g_*[N]_\mathbb{K}=i_*[X]_\mathbb{K}$ in $H_n(X';\mathbb{K})$. Together with the epimorphism $\phi':\pi_1(X')\twoheadrightarrow\pi$ that is induced by the inclusion $X'\hookrightarrow X(n-1)$ this is an extension of $(X,\phi)$. We finally obtain $h:=\hat g:N\to X'$ having the asserted properties.
\end{proof}

\begin{rem*}
Theorem \ref{exists_deg_1_map} does not hold in dimension two. For example consider a closed oriented surface $\Sigma$ of genus $g\geq 2$ and the torus $T^2$. Let $\Phi:\Sigma\to T^2$ be a degree one map. Since the torus is aspherical we are in the situation of case (i) of the theorem. Let $X$ be an extension of $\Sigma$, i.\,e.\ it is obtained by attaching finitely many $1$-cells to $\Sigma$. It is easy to see that there is no map $h:T^2\to X$ that induces a non-trivial homomorphism in $2$-dimensional homology.
\end{rem*}

\begin{proof}[Proof of Theorem \ref{main_thm}]
By the weak comparison axiom and the extension axiom we find
\[ I_\psi(N) \leq I_{\phi'}(X') = I_\phi(X). \]
In case $X=M$ is a manifold we get the equalities of (i) by changing the roles of $M$ and $N$.
\end{proof}


\section{Applications of homological invariance}\label{appl}

Again let $I$ be an invariant that fulfills the weak comparison axiom and the extension axiom, e.\,g.\ $I\in\{\sigma,\sigma^{st},\lambda^n,T\}$. In this section we will apply homological invariance in different situations. First we will demonstrate that orientation-true degree one maps preserve the values of those invariants. This simplifies a rather long proof in \cite{KKM(2006)}. Furthermore, we will prove Corollary \ref{fund_grp_two_elem_intro} about manifolds whose fundamental group consists of only two elements.

\subsection{Degree one maps}

Theorem \ref{main_thm} has the following immediate consequence which improves the homotopy invariance from Corollary \ref{homotopy_invar_2}. 

\begin{defn}
A map $f:M\to N$ between manifolds is called \emph{orientation-true} if it maps orientation preserving loops to orientation preserving ones and orientation reversing loops to orientation reversing ones.
\end{defn}

\begin{cor}\label{deg_one_invar}
If $f:M\to N$ is an orientation-true map of absolute degree one between two connected closed manifolds of dimension $n\geq 3$, then
\[ I_{\psi\circ f_*}(M) = I_\psi(N) \]
for any epimorphism $\psi:\pi_1(N)\twoheadrightarrow\pi$.
\end{cor}

Recall that maps of absolute degree one are always surjective on fundamental groups (see \cite{Epstein(1966)}, Corollary 3.4). Using this, Corollary \ref{deg_one_invar} is a direct application of Theorem \ref{main_thm}.

In \cite{KKM(2006)}, K\c edra, Kotschick, and Morita proved the following theorem (Theorem 4):

\begin{thm*}[K\c edra, Kotschick, Morita]
Let $M$ be a closed oriented manifold with non-vanishing volume flux group $\Gamma_\mu$. Then $M$ has a finite covering $\bar M$ whose minimal volume entropy $\lambda(\bar M)$ vanishes.
\end{thm*}

Their proof on pages 1260--1264 may be shortened and simplified in the following way: it starts with the construction of a map $\Phi:S^1\times F\to \bar M$ from a closed oriented product manifold $S^1\times F$ to a finite covering $\bar M$ of $M$. Lemma 24 on page 1261 states that $\Phi$ has degree one and induces an isomorphism on fundamental groups. Therefore $\lambda(\bar M)=\lambda(S^1\times F)$ by Corollary \ref{deg_one_invar}. Since the minimal volume entropy of $S^1\times F$ vanishes by the vanishing of the minimal volume (see \cite{KKM(2006)} for details), this proves the cited theorem.

\subsection{Adding simply-connected summands}

It is rather difficult to investigate the behaviour of the invariants under connected sums. The easiest case is when one of the summands is simply connected. The next corollary was already known for $I=\lambda$ and $I=\sigma$ by \cite{Babenko(1995)}.

\begin{cor}
Let $M$ and $N$ be two connected closed manifolds and let $\phi:\pi_1(M)\twoheadrightarrow\pi$ be an epimorphism. If $N$ is simply connected, then
\[ I_\phi (M\# N) = I_\phi (M). \]
\end{cor}

\begin{proof}
Let $n:=\dim M=\dim N$. Note that for $n=2$ we necessarily have $N\cong S^2$ and thus $M\# N\cong M$. So we may assume $n\geq 3$. Since the map $M\# N\to M$ is orientation-true, has absolute degree one, and induces an isomorphism on the fundamental group the respective values of $I$ are equal by Corollary \ref{deg_one_invar}.
\end{proof}

\subsection{$\mathbb{Z}_2$-systoles}

Next, we want to look at the case $\phi:\pi_1(M)\twoheadrightarrow\mathbb{Z}_2$ and in particular at manifolds with fundamental group $\mathbb{Z}_2$. Since for finite Galois groups the stable systolic constant, the minimal volume entropy, and the spherical volume vanish, we concentrate on the systolic constant in this subsection. Thus we may use the classification from Corollary \ref{systolic_classification}. Denote $\sigma_n := \sigma(\mathbb{R}P^n)$.

\begin{cor}
Let $M$ be a connected closed manifold of dimension $n\geq 3$ and let $\phi:\pi_1(M)\twoheadrightarrow\mathbb{Z}_2$ be an epimorphism. Then $\sigma_\phi(M)\leq \sigma_n$ and $\sigma_\phi(M)=0$ if and only if $\Phi_*[M]_{\mathbb{Z}_2}=0$. Moreover, if $\tilde M_\phi$ is orientable, then 
\[ \sigma_\phi(M)=
\begin{cases}
\sigma_n & \Phi_*[M]_{\mathbb{Z}_2}\neq 0 \\
0 & \Phi_*[M]_{\mathbb{Z}_2}=0
\end{cases}.\]
In particular, $\sigma_\phi(M)=0$ for $M$ orientable and $n$ even, and also for $M$ non-orientable, $\tilde M_\phi$ orientable and $n$ odd.
\end{cor}

\begin{proof}
Note that $K(\mathbb{Z}_2,1)=\mathbb{R}P^\infty$. Hence
\[ H_n(\mathbb{Z}_2;\mathbb{Z}) = 
\begin{cases}
\mathbb{Z}_2 & n\;\mbox{odd} \\
0 & n>0\;\mbox{even}
\end{cases}\]
and
\[ H_n(\mathbb{Z}_2;\mathcal{O}_{\Id}) = 
\begin{cases}
0 & n\;\mbox{odd} \\
\mathbb{Z}_2 & n>0\;\mbox{even}
\end{cases}.\]
Furthermore $H_n(\mathbb{Z}_2;\mathbb{Z}_2)=\mathbb{Z}_2$ for all $n\geq 0$.

Note that $\Phi_*[M]_\mathbb{K}=0$ if and only if $\Phi_*[M]_{\mathbb{Z}_2}=0$. Hence $\sigma_\phi(M)=0$ if and only if $\Phi_*[M]_{\mathbb{Z}_2}=0$ by the classification of $\phi$-systolic manifolds. Moreover $\Phi_*[M]_\mathbb{K}=0$ in the two particular cases mentioned at the end of the corollary.

If $\Phi_*[M]_{\mathbb{Z}_2}\neq 0$, then $\Phi_*[M]_\mathbb{K}=i_*[\mathbb{R}P^n]_\mathbb{K}$. Hence Theorem \ref{main_thm} finishes the proof.
\end{proof}

Apart from the statement that $\sigma_\phi(M)\in\{0,\sigma_n\}$ for $M$ non-orientable and $\tilde M_\phi$ orientable, this was already proved by Babenko in \cite{Babenko(2004)}. In the special case $\pi_1(M)=\mathbb{Z}_2$ we get Corollary \ref{fund_grp_two_elem_intro}. Note that this statement also holds for $n=2$ since here $M\cong \mathbb{R}P^2$.


\section{Strong extension axiom}\label{str_ext_ax_section}

In the extension axiom of section \ref{ext_ax_section} it was assumed that the simplicial complex $X$ is given together with a surjective homomorphism $\phi:\pi_1(X)\twoheadrightarrow\pi$. This guarantees that the Galois groups of the coverings $\tilde X'_{\phi'}$ and $\tilde X_\phi$ coincide. Without this assumption the Galois group of the extended complex can become extremely large compared to the original one. Think for example of $X=T^n$ with $\phi:\mathbb{Z}^n\hookrightarrow\mathbb{Z}^n*\mathbb{Z}$ the inclusion in the first factor and $X'=T^n\vee S^1$ with $\phi'$ the identity. Here the Galois group of $\tilde X_\phi$ has polynomial growth whereas the one of $\tilde X'_{\phi'}$ grows exponentially. Nevertheless we can show that the systolic constant and the minimal volume entropy behave well in this situation.

\begin{str_ext_ax*} 
Let $(X',\phi')$ be an extension of $(X,\phi)$. Then
\[ I_{\phi'}(X') = I_\phi(X). \]
\end{str_ext_ax*}

We will show that the systolic constant and the minimal volume entropy fulfill this axiom. This is quite easy for the systolic constant, but in case $I=\lambda$ there is some effort necessary. We will approximate Riemannian norms $L_{g,x}$ (see the remark after Definition \ref{norms}) by simpler and more regular norms. This idea is due to Manning (see \cite{Manning(2005)}).

\begin{defn}\label{generator_norms}
Let $G$ be a finitely generated group and $S\subset G$ be finite generating set that is symmetric, i.\,e.\ $S^{-1}=S$. In the following all generating sets will be assumed symmetric. For a norm $L:G\to [0,\infty)$ one defines another norm
\[ N_{L,S} (g) := \inf \{ {\textstyle\sum_{i=1}^n} L(s_i) | g = s_1\cdots s_n, s_i \in S \}. \]
These norms will be called \emph{generator norms}.
\end{defn}

\begin{rem*}
If one takes $L=1$ the norm that assigns $1$ to each non-trivial element of $G$, then $N_{1,S}$ is the well-known \emph{word norm} or \emph{word length} on $G$ with respect to $S$.
\end{rem*}

\begin{lem}\label{entr_inf}
The volume entropy (see Definition \ref{norms}) of generator norms is well-defined and we have
\[ \lambda(G,N_{L,S}) = \inf_{t>0} ({\textstyle\frac{1}{t}}\log\beta_{N_{L,S}}(t) + {\textstyle\frac{1}{t}}\log \# S). \]
\end{lem}

\begin{proof}
Write $\beta := \beta_{N_{L,S}}$. We have
\[ \beta(r+t) \leq \beta(r)\beta(t)\# S .\]
Namely if $N_{L,S}(g)\leq r+t$, then choose a minimal representation $g=s_1\cdots s_n$, i.\,e.\ one fulfilling $N_{L,S}(g)=\sum_{i=1}^n L(s_i)$. Choose $k\in\mathbb{N}_0$ such that
\[ \sum_{i=1}^k L(s_i) \leq r \;\;\;\mbox{and}\;\;\; \sum_{i=1}^{k+1} L(s_i) > r. \]
Define $g_1 := s_1\cdots s_k$ and $g_2 := s_{k+2}\cdots s_n$. Then $g=g_1 s_{k+1} g_2$ and we have
\begin{align*}
N_{L,S}(g_1) & = \sum_{i=1}^k L(s_i) \leq r \;\;\;\mbox{and} \\
N_{L,S}(g_2) & = \sum_{i=1}^n L(s_i) - \sum_{i=1}^{k+1} L(s_i) \leq r+t-r=t.
\end{align*}

Now let $r$ and $t$ be arbitrary positive real numbers and choose $k\in\mathbb{N}_0$ such that $kt<r\leq (k+1)t$. Then
\[ \beta(r) \leq \beta((k+1)t) \leq \beta(kt)\beta(t)\# S \leq\ldots\leq \beta(t)\beta(t)^k\# S^k \]
and consequently
\[ {\textstyle\frac{1}{r}}\log\beta(r) \leq {\textstyle\frac{1}{r}}\log\beta(t) + {\textstyle\frac{k}{r}}\log(\beta(t)\# S) \leq {\textstyle\frac{1}{r}}\log\beta(t) + {\textstyle\frac{1}{t}}\log(\beta(t)\# S). \]
Therefore 
\[ \limsup_{r\to\infty} \textstyle{\frac{1}{r}}\log\beta(r)\leq {\textstyle\frac{1}{t}}\log(\beta(t)\# S) \] 
for all $t>0$. Hence
\[ \limsup_{r\to\infty} {\textstyle\frac{1}{r}}\log\beta(r)\leq \inf_{t>0} ({\textstyle\frac{1}{t}}\log\beta(t) + {\textstyle\frac{1}{t}}\log\# S) \leq \liminf_{t\to\infty} {\textstyle\frac{1}{t}}\log\beta(t). \]
Thus $\lambda(G,N_{L,S})=\lim_{r\to\infty} \frac{1}{r}\log\beta(r)$ exists and fulfills the claimed equality.
\end{proof}

To prove the strong extension axiom for the minimal volume entropy we have to consider the case of one attached circle. The idea is to let its length grow to infinity. The following proposition investigates the analogous situation for generator norms.

\begin{prop}\label{lim_free_prod_gen}
Let $G$ and $H$ be finitely generated groups and $L_G$ and $L_H$ generator norms with respect to the finite generating sets $S$ and $T$, i.\,e.\ $L_G=N_{L_G,S}$ respectively $L_H=N_{L_H,T}$. Then $L_G * \varrho L_H$ is a generator norm on $G*H$ with respect to $S\cup T$ for every $\varrho>0$. We have
\[ \lim_{\varrho\to\infty} \lambda (G*H,L_G*\varrho L_H) = \lambda(G,L_G). \]
\end{prop}

\begin{proof}
Let $\varepsilon>0$. Choose $R>0$ such that
\begin{gather*}
{\textstyle\frac{1}{R}}\log\beta_{L_G}(R) \leq\lambda(G,L_G) +\varepsilon \;\;\;\mbox{and} \\
{\textstyle\frac{1}{R}}\log(\#S+\#T) \leq \varepsilon.
\end{gather*}
For every $\varrho >R/\min_{t\in T}L_H(t)$ we have
\[ \beta_{L_G*\varrho L_H}(R) = \beta_{L_G}(R) \]
since no elements of $H\setminus 1$ are involved, yet. Hence by Lemma \ref{entr_inf}
\begin{align*}
\lambda(G*H,L_G*\varrho L_H) &\leq {\textstyle\frac{1}{R}}\log\beta_{L_G*\varrho L_H}(R)+ {\textstyle\frac{1}{R}}\log(\#S+\#T) \\
&\leq \lambda(G,L_G) +2\varepsilon.
\end{align*}
Therefore $\limsup_{\varrho\to\infty}\lambda(G*H,L_G*\varrho L_H) \leq \lambda(G,L_G)$. But 
\[ \lambda(G*H,L_G*\varrho L_H) \geq \lambda(G,L_G) \] 
is obvious. Thus the limit exists and equals $\lambda(G,L_G)$.
\end{proof}

In the next proposition we want to use Manning's approximation result from \cite{Manning(2005)}. Since we want to swap the limit $\varrho\to\infty$ with the approximation, we have to control the quality of the approximation.

\begin{prop}\label{main_geom-grp_prop}
Let $(X,g)$ be a connected finite Riemannian simplicial complex and $\phi:\pi_1(X)\to\pi$ a homomorphism. As usual denote by $L_{g,x}$ the induced norm on the Galois group $\Gamma:=\pi_1(X)/\ker\phi$. Furthermore let $L_H$ be a generator norm on a finitely generated group $H$ with respect to the finite generating set $T\subset H$. Then
\[ \lim_{\varrho\to\infty} \lambda (\Gamma *H,L_{g,x}*\varrho L_H) = \lambda(\Gamma,L_{g,x}). \]
\end{prop}

\begin{proof}
Write $L_g := L_{g,x}$. Choose a fundamental domain $F\subset\tilde X_\phi$ with diameter $D$ and $x\in F$. Let $R$ be an arbitrary positive real number. Write
\[ h=\gamma_0 h_1\gamma_1\cdots h_n\gamma_n h_{n+1} \in\Gamma * H \]
with $\gamma_0\in\Gamma,\gamma_1,\ldots,\gamma_n\in\Gamma\setminus 1,h_1,\ldots,h_n\in H\setminus 1,h_{n+1}\in H$ and choose $k\in\mathbb{N}_0$ such that
\[ (k-1)R < L_g*\varrho L_H (h) = \sum_{i=0}^n L_g (\gamma_i) + \varrho\cdot\sum_{i=1}^{n+1} L_H(h_i) \leq kR. \]
Now let $k_0,\ldots,k_n\in\mathbb{N}_0$ such that
\[ (k_i -1)R<L_g(\gamma_i) \leq k_i R. \]
Think of the $\gamma_i$ as paths in $\tilde X_\phi$ starting at $x$. Pick points $\alpha_{ij}\in\Gamma,j=1,\ldots,k_i-1$ such that
\[ d(\gamma_i(jR),\alpha_{ij}x) \leq D \]
and set $\alpha_{i0}=1,\alpha_{ik_i}=\gamma_i$. Then
\[ L_g(\alpha_{ij}^{-1}\alpha_{i,j+1}) = d(\alpha_{ij}x,\alpha_{i,j+1}x) \leq R+2D. \]
Put $S:=\{ \alpha\in\Gamma | L_g(\alpha)\leq R+2D \}$. (This is a finite generating system of $\Gamma$ as we just have shown.) Then
\[ N_{L_g,S} (\gamma_i) \leq k_i\cdot (R+2D) \]
and
\begin{align*}
N_{L_g,S}*\varrho L_H (h)
&\leq \sum_{i=0}^n k_i(R+2D) + \varrho\cdot\sum_{i=1}^{n+1}L_H(h_i) \\
&\leq \sum_{i=0}^n (k_i(R+2D)-(k_i-1)R)+ \sum_{i=0}^n L_g (\gamma_i) + \varrho\cdot\sum_{i=1}^{n+1} L_H(h_i) \\
&\leq \sum_{i=0}^n (k_i 2D+R) +kR \\
&\leq (k+n+1)(R+2D)
\end{align*}
since $\sum_{i=0}^n k_i\leq k+n+1$. Hence 
\[ \beta_{L_g*\varrho L_H} (kR) \leq \beta_{N_{L_g,S}*\varrho L_H}((k+n+1)(R+2D)) \] 
and thus
\begin{align*}
\limsup_{r\to\infty} {\textstyle\frac{1}{r}}\log\beta_{L_g*\varrho L_H} (r) 
&= \limsup_{k\to\infty} {\textstyle\frac{1}{kR}}\log\beta_{L_g*\varrho L_H} (kR) \\
&\leq \lim_{k\to\infty} {\textstyle\frac{1}{kR}}\log\beta_{N_{L_g,S}*\varrho L_H} ((k+n+1)(R+2D)) \\
&= {\textstyle\frac{R+2D}{R}} \lambda (\Gamma*H,N_{L_g,S}*\varrho L_H).
\end{align*}
Since for any norm $L\leq N_{L,S}$ by the triangle inequality, we see that 
\[ \liminf_{r\to\infty} {\textstyle\frac{1}{r}}\log\beta_{L_g*\varrho L_H} (r)\geq \lambda(\Gamma*H,N_{L_g,S}*\varrho L_H). \] 
Thus the volume entropy of $L_g*\varrho L_H$ exists and equals $\sup_S \lambda(\Gamma*H,N_{L_g,S}*\varrho L_H)$.

Moreover,
\[ \lambda(\Gamma *H,L_g *\varrho L_H)\leq {\textstyle\frac{R+2D}{R}} \lambda (\Gamma*H,N_{L_g,S}*\varrho L_H). \]
Now, if $\varrho\to\infty$, then by Proposition \ref{lim_free_prod_gen} the right-hand side goes to
\[ {\textstyle\frac{R+2D}{R}} \lambda (\Gamma,N_{L_g,S}) \leq {\textstyle\frac{R+2D}{R}} \lambda (\Gamma,L_g). \]
Since $R>0$ was arbitrary, we get
\[ \limsup_{\varrho\to\infty} \lambda (\Gamma *H,L_g*\varrho L_H) \leq \lambda(\Gamma,L_g). \]
But $\lambda (\Gamma *H,L_g*\varrho L_H) \geq \lambda(\Gamma,L_g)$ is again obvious. Hence the limit exists and equals $\lambda(\Gamma,L_g)$.
\end{proof}

\begin{thm}
The strong extension axiom is fulfilled by the systolic constant $I=\sigma$ and by the minimal volume entropy $I=\lambda$.
\end{thm}

\begin{proof}
By the comparison axiom we have $I_\phi(X) \leq I_{\phi'}(X')$ since the inclusion $X\hookrightarrow X'$ is $(n,1)$-monotone. We may proceed by induction, attaching one cell at a time. The case where the Galois group $\Gamma:=\pi_1(X)/\ker\phi$ does not change is already covered by the extension axiom since here $\phi$ and $\phi'$ factor into epimorphisms onto this quotient and the induced inclusion $\Gamma\hookrightarrow\pi$.

The case remaining to be investigated is therefore the following: Consider $X'=X\cup_h D^1$ where $h:S^0\to X$ is simplicial together with an extension $\phi':\pi_1(X')\cong \pi_1(X)*\mathbb{Z}\to\pi$ of the given homomorphism $\phi$.

First consider the simplicial complex $Y:=X\vee S^1$ where the circle is attached at $h(1)\in X$. Using a path from $h(-1)$ to $h(1)$ in $X$, we get a homotopy equivalence $f:X'\to Y$ and thus may define $\psi:=\phi'\circ f_*^{-1}$. Notice that $(Y,\psi)$ is again an extension of $(X,\phi)$ and that $f$ is $(n,1)$-monotone and has a $(n,1)$-monotone homotopy inverse. Hence $I_{\phi'}(X') = I_\psi(Y)$ and it remains to show that $I_\psi(Y)\leq I_\phi(X)$.

\emph{Case 1:} $I=\sigma$. Let $g$ be a Riemannian metric on $X$. Extend it over $Y$ by assigning the length $\sys_\phi(X,g)$ to the attached circle. Then both $\sys_\psi(Y,g)=\sys_\phi(X,g)$ and $\Vol(Y,g)=\Vol(X,g)$ and consequently $\sigma_\psi(Y)\leq\sigma_\phi(X)$. (If the $\phi$-systole of $(X,g)$ is not finite, i.\,e.\ if $\ker\phi=\pi_1(X)$, then use a sequence of metrics where the length of $S^1$ tends to infinity.)

\emph{Case 2:} $I=\lambda$. Again let a Riemannian metric $g$ on $X$ be given. Define $g_\varrho$ to be the extension over $Y$ that assigns the length $\varrho>0$ to the circle $S^1$. Then $\Vol(Y,g_\varrho)=\Vol(X,g)$. 

Take the attaching point $x=h(1)$ as base point. Then $\pi_1(Y)=\pi_1(X)*\mathbb{Z}$. Consider the homomorphism
\[ \pi_1(Y)\xrightarrow{\phi *\Id} \pi*\mathbb{Z}. \]
We see that $\ker(\phi*\Id)\subset\ker\psi$ and thus
\[ \lambda_{\phi*\Id}(Y) \geq \lambda_\psi(Y). \]
With $L:\mathbb{Z}\to[0,\infty)$ denoting the standard word norm we have
\begin{align*}
\lambda_{\phi*\Id}(Y,g_\varrho) \;\;\;=\;\;\;& \lambda(\Gamma *\mathbb{Z},L_{g,x} *\varrho L) \\
\xrightarrow{\varrho\to\infty}& \lambda (\Gamma, L_{g,x}) = \lambda_\phi (X,g)
\end{align*}
by Proposition \ref{main_geom-grp_prop}. Thus $\lambda_\phi(X) \geq \lambda_{\phi*\Id}(Y) \geq \lambda_\psi (Y)$.
\end{proof}


\section{An inequality between the systolic constant and the minimal volume entropy}\label{ineq_sys_entr}

There is an inequality linking the spherical volume and the minimal volume entropy:
\[ 2^n n^{n/2} T_\phi(M) \leq \lambda_\phi (M)^n. \] 
This was proved by Besson, Courtois, and Gallot in \cite{BCG(1991)}, Th\'eor\`eme 3.8. In this section we will investigate the relation between the systolic constant and the minimal volume entropy. In doing so we will prove a relative analogue of Theorem \ref{sabourau_ineq_intro}. 

\begin{defn}
Let $\pi$ be a group, $\rho:\pi\to\mathbb{Z}_2$ a homomorphism, and $a\in H_n(\pi;\mathcal{O}_\rho)$. Then define
\[ I(a) := \inf_{(X,\Psi)} I_{\Psi_*}(X), \]
where the infimum is taken over all \emph{geometric cycles} $(X,\Psi)$ representing the homology class $a$, i.\,e.\ over all maps $\Psi:X\to K(\pi,1)$ from a connected closed $n$-dimensional pseudomanifold $X$ to $K(\pi,1)$ with $H_n(X;\mathcal{O}_{\rho\Psi_*})\neq 0$ and $\Psi_*[X]_{\mathcal{O}_{\rho\Psi_*}}=a$. For coefficients in $\mathbb{Z}_2$ we use the analogous definition.
\end{defn}

If there is a geometric cycle that is defined on a manifold and that is surjective on fundamental groups, then it is minimal:

\begin{thm}\label{char_min}
If $I$ fulfills the strong extension axiom and the weak comparison axiom, then 
\[ I(\Phi_*[M]_\mathbb{K}) = I_{\Phi_*}(M) \]
for any connected closed manifold of dimension $n \geq 3$ with $\Phi:M\to K(\pi,1)$ such that $\Phi_*:\pi_1(M)\twoheadrightarrow \pi$ is surjective. (We use $\mathbb{K}$ the orientation bundle of $M$ if the covering $\tilde M_{\Phi_*}$ is orientable and $\mathbb{K}=\mathbb{Z}_2$ otherwise.)
\end{thm}

\begin{proof}
Let $(X,\Psi)$ be a geometric cycle representing $\Phi_*[M]_\mathbb{K}$. Then there exists another geometric cycle $(X',\Psi')$ representing the same homology class such that $\Psi'$ maps the fundamental group of $X'$ surjectively onto $\pi$ and such that $I_{\Psi'_*}(X')= I_{\Psi_*}(X)$. 

Namely, let $\gamma_1,\ldots,\gamma_m$ be generators of $\pi$. Think of the $\gamma_i$ as closed curves in $K(\pi,1)$ and define
\begin{gather*} 
X'':= X\vee\left(\bigvee_{i=1}^m S^1\right) \;\;\;\mbox{and}\\ 
\Psi'':=\Psi\vee\left(\bigvee_{i=1}^m \gamma_i\right): X'' \to K(\pi,1).
\end{gather*}
Then $\Psi''_*[X'']_\mathbb{K}=\Phi_*[M]_\mathbb{K}$ (here $[X'']_\mathbb{K}$ is the image of $[X]_\mathbb{K}$ under the inclusion $X\hookrightarrow X''$) and the induced homomorphism $\Psi''_*$ on the fundamental group is an epimorphism. By the strong extension axiom $I_{\Psi''_*}(X'') = I_{\Psi_*}(X)$. 

Now consider the pseudomanifold
\[ X' := X \;\#\; \left(\mathop{\#}\limits_{i=1}^m (S^n/\{\pm \operatorname{pt}\})\right), \]
where $S^n/\{\pm \operatorname{pt}\}$ is the $n$-sphere with two points identified. The projection of $S^n$ to a closed interval such that $\pm \operatorname{pt}$ are mapped to the boundary points induces a map $S^n/\{\pm \operatorname{pt}\}\to S^1$. Let $p:X'\to X''$ be the composition of the projection
\[ X' \to  X\vee\left(\bigvee_{i=1}^m S^n/\{\pm \operatorname{pt}\}\right) \]
with this map on each $S^n/\{\pm \operatorname{pt}\}$. Define 
\[ \Psi':=\Psi''\circ p :X'\to K(\pi,1). \]
Note that $p$ is a homotopy equivalence. Thus $\Psi'$ induces a surjection on fundamental groups, $H_n(X';\mathbb{K})\neq 0$, and $(X',\Psi')$ represents $\Phi_*[M]_\mathbb{K}$. 

Since $p$ can be chosen strictly $(n,1)$-monotone and has a strictly $(n,1)$-monotone homotopy inverse, we get
\[ I_{\Psi'_*}(X') = I_{\Psi''_*}(X'') = I_{\Psi_*}(X) \]
by the weak comparison axiom. From Theorem \ref{exists_deg_1_map} it follows that $I_{\Phi_*}(M) \leq I_{\Psi'_*}(X')= I_{\Psi_*}(X)$.
\end{proof}

Now we follow \cite{Sabourau(2006)}, sections 4 and 5. First we need a theorem of Gromov:

\begin{thm}[\cite{Gromov(1983)}, pages 70, 71]
There exists a constant $A_n >0$ such that for all $\varepsilon>0$ there is a geometric cycle $(X,\Psi)$ representing $a\in H_n(\pi;\mathbb{K})$ and a Riemannian metric $g$ on $X$ such that
\[ \sigma_{\Psi_*}(X,g) \leq (1+\varepsilon) \sigma(a) \]
and 
\[ \Vol B(x,R) \geq A_n R^n \]
for all $x\in X$ and $\varepsilon\leq R/\sys_{\Psi_*}(X,g) \leq \frac{1}{2}$. (Here, $B(x,R)$ denotes the ball around $x\in X$ of radius $R$ with respect to $g$.) Such cycles are called \emph{$\varepsilon$-regular}.
\end{thm}

In fact, Gromov proved this only for $\mathbb{K}=\mathbb{Z}$ and $\mathbb{K}=\mathbb{Z}_2$ but it remains true for local integer coefficients with the same arguments. On $\varepsilon$-regular geometric cycles one may compare the systolic constant and the minimal volume entropy:

\begin{prop}[\cite{Sabourau(2006)}, Proposition 4.1]
Let $(X,g,\Psi)$ be an $\varepsilon$-regular geometric cycle. Then
\[ \lambda_{\Psi_*}(X,g) \Vol(X,g)^{1/n} \leq \frac{\sigma_{\Psi_*}(X,g)^{1/n}}{\beta} \log \frac{\sigma_{\Psi_*}(X,g)}{A_n \alpha^n} \]
for all $\alpha\geq \varepsilon,\beta>0$ with $4\alpha+\beta <\frac{1}{2}$. 
\end{prop}

From this it follows directly that
\[ \lambda(a) \leq \frac{\sigma(a)^{1/n}}{\beta} \log \frac{\sigma(a)}{A_n \alpha^n} \]
for all $\alpha,\beta>0$ with $4\alpha+\beta <\frac{1}{2}$. The calculation from the proof of \cite{Sabourau(2006)}, Theorem 5.1 shows:

\begin{cor}
There exists a constant $c_n>0$ such that
\[ \sigma(a) \geq c_n \frac{\lambda(a)^n}{\log^n(1+\lambda(a))}. \]
\end{cor}

Combined with Theorem \ref{char_min} this proves:

\begin{thm}
Let $M$ be a connected closed manifold of dimension $n\geq 3$ and $\phi:\pi_1(M)\twoheadrightarrow\pi$ an epimorphism. Then there exists a positive constant $c_n$ depending only on $n$ such that
\[ \sigma_\phi (M) \geq c_n \frac{\lambda_\phi (M)^n}{\log^n(1+\lambda_\phi(M))}. \]
\end{thm}

Special cases of this statement were shown by Sabourau (see \cite{Sabourau(2006)}). Note that the absolute version of this theorem is in fact Theorem \ref{sabourau_ineq_intro} for $n\geq 3$. The two-dimensional case was proved in \cite{KS(2005)}. Thus Theorem \ref{sabourau_ineq_intro} is shown.


\providecommand{\bysame}{\leavevmode\hbox to3em{\hrulefill}\thinspace}
\providecommand{\MR}{\relax\ifhmode\unskip\space\fi MR }
\providecommand{\MRhref}[2]{%
  \href{http://www.ams.org/mathscinet-getitem?mr=#1}{#2}
}
\providecommand{\href}[2]{#2}

\end{document}